\documentclass[12pt,oneside]{article}
\usepackage{amsmath,amssymb,amsfonts,amsthm}
\usepackage{color}
\def\goth{\mathfrak}

\newcommand{\T}{{\cal T}}

\newcommand{\prof}{\noindent \textit{\textbf{Proof.\:\:}}}
\newcommand{\tm}{\T M}
\newcommand{\p}{\pi^{-1}(TM)}
\newcommand{\cp}{\mathfrak{X}(\pi (M))}

\def\o#1{\overline{#1}}
\def\ti#1{\tilde{#1}}

\def\x{{\goth X}}

\def\Section#1{\vspace{30truept}\addtocounter{section}{1}\setcounter{thm}{0}\setcounter{equation}{0}
{\noindent\Large\bf\arabic{section}.~~#1}\par \vspace{12pt}}

\newtheorem{thm}{Theorem}[section]

\newtheorem{lem}[thm]{Lemma}
\newtheorem{prop}[thm]{Proposition}

\numberwithin{equation}{section}

\title{TWO NONRELATED FINSLER STRUCTURES ON A MANIFOLD}
\author{Aly A. Tamim\and Nabil L. Youssef}
\date{}

\begin{document}
\bibliographystyle{plain}
\maketitle \vspace{-1.15cm}
\begin{center}
{Department of Mathematics, Faculty of Science,\\ Cairo University,
Giza, Egypt.}
\end{center}
\maketitle

\begin{center}
Dedicated to Professor Radu Miron\\
on the occasion of his 70th birthday
\end{center}

\vspace{20truept}\centerline{\Large\bf{Introduction}}\vspace{12pt}
\par
It is a standard and commonly used approach to consider two
different differentiable structures on the same base manifold, which
are related a priori in some sort (projectively,
conformally,...etc.) and to study the relationships between the
corresponding geometric objects (connection, geodesic, curvature
tensors,...etc.) associated to those structures. In this direction,
the reader may refer for example to~\cite{[1]},~\cite{[3]} and\,
\cite{[7]}. In the present paper, we proceed differently: We
consider two different {\em Finsler} structures $L$ and $L^*$ on the
same base manifold $M$, with no relation preassumed between them.
\par
Introducing the $\pi$-tensor field representing the difference
between the Cartan's connections associated with $L$ and $L^*$, we
investigate the (necessary and sufficient) conditions, to be
satisfied by this $\pi$-tensor field, for the geometric objects
associated with $L$ and $L^*$ to have the same properties. Among
various items investigated in the paper, we consider the properties
of being a geodesic, a Jacobi field, a Berwald manifold, a locally
Minkowskian manifold and a Landsberg manifold.
\par
It should be noticed that our approach is a global one. That is, it
does not make use of local coordinate techniques.


\Section{Notations and Preliminaries}
 In this section we give a
brief account of the basic concepts necessary for this work. For
more details, refer to~\cite{[2]} or ~\cite{[5]}. We make the
general assumption that all geometric objects we consider are of
class $C^\infty$. The following notations will be used throughout
the paper:\newline \noindent $M$: a differentiable manifold of
finite dimension and of class $C^\infty$.\newline
$\pi_M:TM\longrightarrow M$: the tangent bundle of $M$.\newline
$\pi:\tm\longrightarrow M$: the subbundle of nonzero
     vectors tangent to $M$.\newline
$P:\pi^{-1}(TM)\longrightarrow \tm$: the bundle, with base
     space $\tm$, induced by $\pi$ and $TM$.\newline
${\goth F}(M)$: the ${\Bbb R}$-algebra of differentiable functions
     on $M$.\newline
${\goth X}(M)$: the ${\goth F}(M)$-module of vector fields on
$M$.\newline $\cp$: the ${\goth F}(\tm)$-module of differentiable
sections of $\pi^{-1}(TM)$.
\par
Elements of $\cp$ will be called $\pi$-vector fields and will be
denoted by barred letters $\o X$. Tensor fields on $\p$ will be
called $\pi$-tensor fields. The fundamental vector field is the
$\pi$-vector field $\o\eta$ defined by $\o\eta(u)=(u,u)$ for all
$u\in\tm$. The lift to $\p$ of a vector field $X$ on $M$ is the
$\pi$-vector field $\o X$ defined by $\o X(u)=(u,X(\pi(u)))$.
\par
The tangent bundle $T(\tm)$ is related to the vector bundle $\p$ by
the exact sequence:
$$0\longrightarrow\p\stackrel{\gamma}\longrightarrow T(\tm)\stackrel{\rho}
\longrightarrow\p\longrightarrow0,$$ where the vector bundle
morphisms are defined by $\rho=(\pi_{\tm},d\pi)$ and
$\gamma(u,v)=j_u(v)$, where $j_u$ is the natural isomorphism
$j_u:T_{\pi_M(v)}M\longrightarrow T_u(T_{\pi_M(v)}M)$.
\par
Let $\nabla$ be an affine connection (or simply a connection) in the
vector bundle $\p$. We associate to $\nabla$ the map
$$K:\tm\longrightarrow\p:X\longmapsto\nabla_X\o\eta,$$
called the connection map of $\nabla$. A tangent vector $X\in
T_u(\tm)$ is said to be horizontal if $\ K(X)=0$. The connection
$\nabla$ is said to be regular if
$$T_u(\tm)=V_u(\tm)\oplus H_u(\tm)\qquad\forall u\in\tm,$$
where $V_u(\tm)$ and $H_u(\tm)$ are respectively the vertical and
horizontal spaces at $u$. If $M$ is endowed with a regular
connection, we can define a section $\beta$ of the morphism $\rho$
by $\ \beta=(\rho\mid_{H(\tm)})^{-1}$. It is clear that $\
\rho\circ\beta\ $ is the identity map on $\p$ and $\ \beta\circ\rho\
$ is the identity map on $H(\tm)$.
\par
For every $X,Y\in\x(\tm)$, the torsion form ${\bf T}$ and the
curvature transformation ${\bf R}$ of the connection $\nabla$ are
defined by:
$${\bf T}(X,Y)=\nabla_X\rho Y-\nabla_Y\rho X-\rho[X,Y],\quad
{\bf R}(X,Y)=[\nabla_Y,\nabla_X]+\nabla_{[X,Y]}.$$ The horizontal
and mixed torsion tensors, denoted respectively by $S$ and $T$, are
defined, for all $\o X,\o Y\in\cp$, by:
$$S(\o X,\o Y)={\bf T}(\beta\o X,\beta\o Y),\quad
T(\o X,\o Y)={\bf T}(\gamma\o X,\beta\o Y).$$ The horizontal, mixed
and vertical curvature tensors, denoted respectively by $R$, $P$ and
$Q$, are defined, for all $\o X, \o Y,\o Z\in\cp$, by:
$$R(\o X,\o Y)\o Z={\bf R}(\beta \o X,\beta\o Y)\o Z,\ \
P(\o X,\o Y)\o Z={\bf R}(\gamma\o X,\beta\o Y)\o Z,\ \ Q(\o X,\o
Y)\o Z={\bf R}(\gamma\o X,\gamma\o Y)\o Z.$$
\par
If $c:I\longrightarrow M$ is a regular curve in $M$, its canonical
lift to $\tm$ is the curve $\ti c$ defined by $\ti c:t\longmapsto
dc/dt$. The lift of a vector field $X$ along $c$ is the $\pi$-vector
field along $\ti c$ defined by $\ \o X:\ti c(t)\longmapsto(\ti
c(t),X(c(t)))$. In particular, the velocity vector field $d c/dt$
along $c$ is lifted to the $\pi$-vector field $\o
dc/dt=(dc/dt,dc/dt)$ along $\ti c$. Clearly, $\rho(d\ti c/dt)=\o
dc/dt =\o\eta\mid_{\ti c(t)}$. A vector field $X$ along a regular
curve $c$ in $M$ is parallel along $c$ with respect to the
connection $\nabla$ if $\ D\o X/dt=0$, where $D/dt$ is the covariant
derivative operator, associated with $\nabla$, along $\ti c$. A
regular curve $c$ in $M$ is a geodesic if the $\pi$-vector field
$D(\o dc/dt)/dt$ vanishes identically. In this case, the vector
field $d\ti c/dt$ along $\ti c$ is horizontal.

\Section{Connections}
 Let $L$ and $L^*$ be two Finsler structures
on a differentiable manifold $M$, with no relation assumed a priori
between them. Throughout this work, entities of $(M,L^*)$ will be
marked by an asterisk "$^*$".
\par
Let $\nabla$ and $\nabla^*$ be the Cartan's connections associated
respectively with the Finsler manifolds $(M,L)$ and $(M,L^*)$. For
every $X\in\x(\tm)$ and $\o Y\in\cp$, let us write

\begin{equation}\label{eq.1}
{\nabla^*}_X\o Y=\nabla_X\o Y+U(X,\o Y),
\end{equation}
where $U$ is an ${\goth F}(\tm)$-bilinear mapping $\x(\tm)\times\cp
\longrightarrow\cp$ representing the difference between the two
connections $\nabla^*$ and $\nabla$.\newline For every $\o X,\o
Y\in\cp$, we set

\begin{equation}\label{eq.2}
  \left.
      \begin{array}{rclcl}
A(\o X,\o Y) &:=& U(\gamma\o X,\o Y),\qquad B(\o X,\o Y) &:=& U(\beta\o X,\o Y)\\
     N(\o X) &:=& B(\o X,\o\eta),\qquad\qquad\quad\  N_0 &:=& N(\o\eta)
      \end{array}
  \right\}
\end{equation}
As a vector field $X$ on $\tm$ can be represented by

\begin{equation}\label{eq.3}
X=\gamma KX+\beta\rho X,
\end{equation}
it follows from (\ref{eq.2}) and (\ref{eq.3}) that

\begin{equation}\label{eq.4}
U(X,\o Y)=A(K(X),\o Y)+B(\rho X,\o Y),\qquad\forall X\in\x(\tm),\ \o
Y\in\cp.
\end{equation}

 By (\ref{eq.1}) and (\ref{eq.2}), taking the regularity of the
connections $\nabla$ and $\nabla^*$ into account, we get

\begin{lem}\label{le.1}
For every $\pi$-vector field $\o X$, we have\:
$$A(\o X,\o\eta)=0.$$
\end{lem}
\par
Using Lemma \ref{le.1}, equations (\ref{eq.1}), (\ref{eq.2}) and
(\ref{eq.4}) imply
\begin{lem}\label{le.2}
The relation between the connection maps $K$ and $K^*$ is given by\:
$$K^*=K+N\circ\rho.$$
\end{lem}

\par
The following result follows directly from Lemma \ref{le.2}.
\begin{prop}\label{pp.1}
A horizontal vector field with respect to $\nabla$ (resp.
$\nabla^*$) is horizontal with respect $\nabla^*$ (resp. $\nabla$)
if, and only if, $N\circ\rho=0$.
\end{prop}
\par
The proof of the following result is not difficult.

\begin{prop}\label{pp.2}
The relation between $\beta$ and $\beta^*$ is given by\:
$$\beta^*=\beta-\gamma\circ N.$$
\end{prop}
\begin{prop}\label{pp.3}
For every $\o X,\o Y\in\cp$, we have\:\newline {\em\textbf{(a)}}
$T^*(\o X,\o Y)=T(\o X,\o Y)+A(\o X,\o Y)$.\newline
{\em\textbf{(b)}} $T^*(N(\o X),\o Y)-T^*(N(\o Y),\o X)=B(\o X,\o
Y)-B(\o Y,\o X)$.
\end{prop}
\prof One can easily show, for every $X,Y\in\x(\tm)$, that
\begin{equation}\label{eq.5}
{\bf T^*}(X,Y)={\bf T}(X,Y)+U(X,\rho Y)-U(Y,\rho X).
\end{equation}
(a) Setting $X=\rho\o X$ and $Y=\beta\o Y$ in (\ref{eq.5}), using
Proposition \ref{pp.2} and taking the definition of $T$ into
account, (a) follows.\newline (b) Setting $X=\beta\o X$ and
$Y=\beta\o Y$ in (\ref{eq.5}), using Proposition \ref{pp.2} and the
fact that $S^*(\o X,\o Y)=0=S(\o X,\o Y)$, (b) follows.\ \ $\Box$


\Section{Curvature Tensors}
 Let ${\bf R}$ and ${\bf R^*}$ be the
curvature transformations of the connections $\nabla$ and $\nabla^*$
respectively.
\begin{lem}\label{le.3}
For every $\ X,Y\in\x(\tm)$, $\ \o Z\in\cp$, we have\:
$$ {\bf R^*}(X,Y)\o Z={\bf R}(X,Y)\o Z+\Omega(X,Y)\o Z,$$
where \newline
\smallskip
$\Omega(X,Y)\o Z = (\nabla_YB)(\rho X,\o Z)-(\nabla_XB)(\rho Y,\o Z)
                    +(\nabla_YA)(K(X),\o Z)\\
\phantom{.\qquad\qquad\quad}  -(\nabla_XA)(K(Y),\o Z)
                    +A({\bf R}(X,Y)\o\eta,\o Z)-B({\bf T}(X,Y),\o Z)\\
\phantom{.\qquad\qquad\quad}  +U(Y,U(X,\o Z))-U(X,U(Y,\o Z)).$
\end{lem}
\par The following proposition gives some useful technical formulas
which will be used frequently in the sequel.
\begin{prop}\label{pp.4}
For every $\ \o X,\o Y,\o Z\in\cp$, we have\:\newline {\em
\textbf{(a)}} $R^*(\o X,\o Y)\o Z+P^*(N(\o X),\o Y)\o Z-P^*(N(\o
Y),\o X)\o Z
      +Q^*(N(\o X),N(\o Y))\o Z\\
      \phantom{\ \qquad\qquad\qquad}=R(\o X,\o Y)\o Z
      +\Omega(\beta\o X,\beta\o Y)\o Z$,\newline
where\newline $\Omega(\beta\o X,\beta\o Y)\o Z
      = (\nabla_{\beta\o Y}B)(\o X,\o Z)
          -(\nabla_{\beta\o X}B)(\o Y,\o Z)+A(R(\o X,\o Y)\o\eta,\o Z)\\
\phantom{..\qquad\qquad\qquad} +B(\o Y,B(\o X,\o Z))-B(\o X,B(\o
Y,\o Z)).$\newline {\em \textbf{(b)}} $P^*(\o X,\o Y)\o Z+Q^*(\o
X,N(\o Y))\o Z=P(\o X,\o Y)\o Z
      +\Omega(\gamma\o X,\beta\o Y)\o Z$,\newline
where\newline $\Omega(\gamma\o X,\beta\o Y)\o Z
      =-(\nabla_{\gamma\o X}B)(\o Y,\o Z)
          +(\nabla_{\beta\o Y}A)(\o X,\o Z)+A(P(\o X,\o Y)\o\eta,\o Z)\\
\phantom{.\ \qquad\qquad\qquad}-B(T(\o X,\o Y),\o Z)+B(\o Y,A(\o
X,\o Z))
          -A(\o X,B(\o Y,\o Z)).$\newline
{\em\textbf{(c)}} $Q^*(\o X,\o Y)\o Z=Q(\o X,\o Y)\o Z
      +\Omega(\gamma\o X,\gamma\o Y)\o Z$,\newline
where\newline $\Omega(\gamma\o X,\gamma\o Y)\o Z=(\nabla_{\gamma\o
Y}A)(\o X,\o Z)
      -(\nabla_{\gamma\o X}A)(\o Y,\o Z)+A(\o Y,A(\o X,\o Z))
      -A(\o X,A(\o Y,\o Z))$.\newline
\par In particular, if $\o Z=\o\eta$, we get\newline
{\em\textbf{(a)}}$'$ $R^*(\o X,\o Y)\o\eta+P^*(N(\o X),\o Y)\o\eta
         -P^*(N(\o Y),\o X)\o\eta
         =R(\o X,\o Y)\o\eta+\Omega(\beta\o X,\beta\o Y)\o\eta$,\newline
where\newline
         $\Omega(\beta\o X,\beta\o Y)\o\eta=(\nabla_{\beta\o Y}N)(\o X)
         -(\nabla_{\beta\o X}N)(\o Y)
         +B(\o Y,N(\o X))-B(\o X,N(\o Y))$.\newline
{\em\textbf{(b)}}$'$ $P^*(\o X,\o Y)\o\eta=P(\o X,\o Y)\o\eta
         +\Omega(\gamma\o X,\beta\o Y)\o\eta$,\newline
where\newline
         $\Omega(\gamma\o X,\beta\o Y)\o\eta=-(\nabla_{\gamma\o X}N)(\o Y)
         +B(\o Y,\o X)-N(T(\o X,\o Y))-A(\o X,N(\o Y))$.\newline
{\em\textbf{(c)}}$'$ $A(\o X,\o Y)=A(\o Y,\o X)$,\newline that is,
the $\pi$-tensor field $A$ is symmetric.
\end{prop}
\prof \newline (a) follws from Lemma \ref{le.3} for $X=\beta\o X,\
Y=\beta\o Y$ and from Proposition \ref{pp.2}.\newline (b) follws
from Lemma \ref{le.3} for $X=\gamma\o X,\ Y=\beta\o Y$ and from
Proposition \ref{pp.2}.\newline (c) follws from Lemma \ref{le.3} for
$X=\gamma\o X,\ Y=\gamma\o Y$.\ \ $\Box$
\begin{prop}\label{pp.5}
Assume that the $\pi$-tensor field $B$ vanishes. Then, $R^*=0$ if,
and only if, $R=0$.
\end{prop}
 \prof Since $B=0$, it follows from Proposition \ref{pp.4}(a) that
$$R^*(\o X,\o Y)\o Z=R(\o X,\o Y)\o Z+A(R(\o X,\o Y)\o\eta,\o Z).$$
It is thus clear that if $R=0$, then $R^*=0$.
\par Conversely, if $R^*=0$, then
\begin{equation}\label{eq.6}
R(\o X,\o Y)\o Z=-A(R(\o X,\o Y)\o\eta,\o Z).
\end{equation}
Setting $\o Z=\o\eta$ in (\ref{eq.6}) and using Lemma \ref{le.1}, we
get $R(\o X,\o Y)\o\eta=0$. Hence, it follows again from
(\ref{eq.6}) that $R=0$.\ \ $\Box$
\begin{prop}\label{pp.6}
The $\pi$-tensor field $N$ vanishes if, and only if, $N_0$ vanishes.
\end{prop}
\prof  Firstly, it is clear that $N=0$ implies $N_0=0$.
\newline Setting $\o Y=\o\eta$ in Proposition \ref{pp.4}(b)$'$ and using the
properties of the torsion and curvature tensors, we get
$$\nabla_{\gamma\o X}N_0=N(\o X)+B(\o\eta,\o X)-A(\o X,N_0). $$
On the other hand, we have from Proposition \ref{pp.3}(b)
$$T^*(N_0,\o X)=B(\o\eta,\o X)-N(\o X).$$
It follows, from the above two identities, that
$$\nabla_{\gamma\o X}N_0=2N(\o X)+T^*(N_0,\o X)-A(\o X,N_0). $$
Consequently, $N_0=0$ implies $N=0$.\ \ $\Box$
\medskip
\par
Now, let us assume that the $\pi$-vector field $N_0$ vanishes. Then,
by Proposition \ref{pp.6}, the formula (a)$'$ of Proposition
\ref{pp.4} takes the form $\ R^*(\o X,\o Y)\o\eta=R(\o X,\o
Y)\o\eta$. It is well-known that~\cite{[6]} the horizontal
distribution with respect to $\nabla$ is completely integrable if,
and only if, $R(\o X,\o Y)\o\eta=0$. Therefore, we have
\begin{thm}\label{th.1}
Suppose that the $\pi$-vector field $N_0$ vanishes. The horizontal
distribution with respect to $\nabla^*$ is completely integrable if,
and only if, the horizontal distribution with respect to $\nabla$ is
completely integrable.
\end{thm}

\Section{Geodesics and Jacobi Fields}
 Let $D/dt$ and $D^*/dt$ be
the covariant derivative operators, corresponding respectively to
the Cartan's connections $\nabla$ and $\nabla^*$, along a curve $\ti
c$ in $\tm$. One can easily show that
\begin{equation}\label{eq.7}
D^*\o X/dt=D\o X/dt+U(d\ti c/dt,\o X),
\end{equation}
for every $\pi$-vector field $\o X$ along $\ti c$. This formula
gives directly
\begin{lem}\label{le.4}
Let $\ti c$ be a curve in $\tm$. A parallel $\pi$-vector field $\o
X$ along $\ \ti c\ $ in $(M,L)$ {\em(resp. $(M,L^*)$)} is parallel
along $\ti c$ in $(M,L^*)$ {\em(resp. $(M,L)$)} if, and only
if,\newline $U(d\ti c/dt,\o X)=0$.
\end{lem}
\begin{thm}\label{th.2}
A necessary and sufficient condition for a geodesic $c$ in $(M,L)$
{\em(resp. $(M,L^*)$)} to be a geodesic in $(M,L^*)$ {\em(resp.
$(M,L)$)} is that $B(\o V,\o V)=0$, where $\o V=\o\eta\mid_{\ti
c(t)}$.
\par In other words, the Finsler manifolds $(M,L)$ and $(M,L^*)$
are projectively related if, and only if, $B(\o V,\o V)=0$ for every
$\o V$.
\end{thm}
\prof  Let $c$ be a regular curve in $M$. The canonical lift $\ti c$
of $c$ to $\tm$ is such that $\rho(d\ti c/dt)=\o\eta\mid_{\ti
c(t)}=\o V$. Then, by (\ref{eq.4}) and Lemma \ref{le.1}, we have $\
U(d\ti c/dt,\o V)=B(\o V,\o V)$. It follows then from (\ref{eq.7})
that ${D^*\o V}/dt=D\o V/dt+B(\o V,\o V)$. The result follows from
this relation and the fact that $c$ is a geodesic in $(M,L)$ (resp.
$(M,L^*)$) if, and only if, $D\o V/dt=0$ (resp. $D^*\o V/dt=0$).\ \
$\Box$
\medskip
\par
A vector field $J$ along a geodesic $c$ in $M$ is called a Jacobi
field with respect to $\nabla$ if it satisfies the Jacobi
differential equation
$$\frac{D^2\o J}{dt^2}+R(\o V,\o J)\o V=0,$$
where $\o J$ and $\o V$ are respectively the lifts of $J$ and
$V=dc/dt$ along $\ti c$.
\par
Writing equation (\ref{eq.7}) for $\o X=\o J$ and using
(\ref{eq.4}), we get
\begin{equation}\label{eq.8}
\frac{D^*\o J}{dt}=\frac{D\o J}{dt}+B(\o V,\o J) +A(K({d\ti
c}/{dt}),\o J).
\end{equation}
Proposition 3.2(a) for $\o X=\o Z=\o V$ and $\o Y=\o J$ yields
\begin{equation}\label{eq.9}
R^*(\o V,\o J)\o V+P^*(N(\o V),\o J)\o V=R(\o V,\o J)\o V
      +\Omega(\beta\o V,\beta\o J)\o V,
\end{equation}
where\newline $\Omega(\beta\o V,\beta\o J)\o V=\nabla_{\beta\o
J}B(\o V,\o V) -(\nabla_{\beta\o V}B)(\o J,\o V) +B(\o J,B(\o V,\o
V))-B(\o V,B(\o J,\o V))$.\newline Now, if $B(\o V,\o V)=0$, it
follows from Proposition \ref{pp.6} that $B(\o X,\o V)=0$ for every
$\pi$-vector field $\o X$. Then, by Proposition \ref{pp.3}, we have
$B(\o V,\o X)=B(\o X,\o V)=0$, and consequently $\Omega(\beta\o
V,\beta\o J)\o V=0$. Moreover, $K({d\ti c}/{dt})=0$ since the vector
field ${d\ti c}/{dt}$ is horizontal with respect to $\nabla$.
Therefore, equations (\ref{eq.8}) and (\ref{eq.9}) reduce
respectively to
$$\frac{D^*\o J}{dt}=\frac{D\o J}{dt}\quad\mbox{ and }\quad
R^*(\o V,\o J)\o V=R(\o V,\o J)\o V.$$ These two identities imply:
$$\frac{{D^*}^2\o J}{dt^2}+R^*(\o V,\o J)\o V=\frac{D^2\o J}{dt^2}+R(\o V,\o J)\o V.$$
Hence, by the Jacobi equation, we obtain
\begin{thm}\label{th.3}
If $B(\o V,\o V)=0$ for every geodesic $c$ in $M$, then the two
Finsler manifolds $(M,L)$ and $(M,L^*)$ have the same Jacobi fields.
\end{thm}
\par
It is not difficult to show that if $B(\o V,\o V)=0$ for every
geodesic $c$ in $M$, then the geodesics of the two Finsler manifolds
$(M,L)$ and $(M,L^*)$ possess the same conjugate points. Therefore,
as a consequence of Theorem \ref{th.3} and the well-known Morse
index theorem [2], we get
\begin{thm}\label{th.4}
If $B(\o V,\o V)=0$ for every geodesic $c$ in $M$, then the
geodesics of the two Finsler manifolds $(M,L)$ and $(M,L^*)$ have
the same Morse index.
\end{thm}

\Section{Special Finsler Manifolds}
 A Finsler manifold $(M,L)$ is a
Berwald manifold~\cite{[3]} if the torsion tensor $T$ satisfies the
condition that $\ \nabla_{\beta\o X}T=0\ $ for every $\o X\in\cp$.
\medskip
\par
Suppose now that the $\pi$-tensor field $B$ vanishes. Then, by
Proposition \ref{pp.2}, we have $\beta^*=\beta$. Hence,
${\nabla^*}_{\beta^*\o X}\o Y=\nabla_{\beta\o X}\o Y\quad\forall \o
X,\o Y\in\cp$. From which, using Proposition \ref{pp.3}, we get
$${\nabla^*}_{\beta^*\o X}T^*=\nabla_{\beta\o X}T+\nabla_{\beta\o X}A.$$
This relation implies
\begin{thm}\label{th.5}
Assume that the $\pi$-tensor field $B$ vanishes. Let $(M,L)$
{\em(resp.\,\,$(M,L^*)$)} be a Berwald manifold. A necessary and
sufficient condition for $(M,L^*)$ {\em(resp. $(M,L)$)} to be a
Berwald manifold is that $\ \nabla_{\beta\o X}A=0\ $ for all $\o
X\in\cp$.
\end{thm}
\smallskip
\par
A Finsler manifold $(M,L)$ is locally Minkowskian~\cite{[3]} if, and
only if, $R=0$ and $\nabla_{\beta\o X}T=0\ $ for all $\o X\in\cp$.
\smallskip
\par
Combining Proposition \ref{pp.5} and Theorem \ref{th.5}, we get
\begin{thm}\label{th.6}
Assume that the $\pi$-tensor field $B$ vanishes. Let $(M,L)$
{\em(resp.$(M,L^*)$)} be a locally Minkowskian manifold. A necessary
and sufficient condition for $(M,L^*)$ {\em(resp. $(M,L)$)} to be
locally Minkowskian is that $\ \nabla_{\beta\o X}A=0\ $ for all $\o
X\in\cp$.
\end{thm}
\smallskip
\par
A Finsler manifold $(M,L)$ is a Landsberg manifold~\cite{[4]} if it
satisfies the condition that $\ P(\o X,\o Y)\o\eta=0\ $ for all $\o
X,\o Y\in\cp$.
\smallskip
\par
Using Proposition \ref{pp.4}(b)$'$, we obtain
\begin{thm}\label{th.7}
Suppose that the $\pi$-tensor field $B$ vanishes. The Finsler
manifold $(M,L)$ is a Landsberg manifold if, and only if, the
Finsler manifold $(M,L^*)$ is a Landsberg manifold.
\end{thm}

\bigskip
\noindent {\bf Concluding Remark.} The results obtained in this
paper may be applied to treat any two related Finsler structures,
where the $\pi$-tensor fields $A$ and $B$ take special forms
depending on the type of relation between the two structures. In a
forthcoming paper, we shall apply these results to "Generalized
Randers Manifolds".
\newpage


\end{document}